%% file: stau.tex
\documentclass[12pt]{amsart}
\font\emailfont=cmtt10

\usepackage{amsmath,amsthm,amsfonts,amscd,flafter,epsf}
\usepackage{graphics}
\usepackage{epsfig}
\usepackage{psfrag}

\input{mac}

\commentable{

\title[{The Ozsv{\'a}th-Szab{\'o} and Rasmussen concordance invariants are not equal}] 
{The Ozsv{\'a}th-Szab{\'o} and Rasmussen concordance invariants are not equal}

}

\author[Matthew Hedden]{Matthew Hedden}
\address{Department of
Mathematics, Princeton University, Princeton, NJ 08540 \newline
\indent{\emailfont{mhedden@math.princeton.edu}}}

\author[Philip Ording]{Philip Ording}
\address{Department of
Mathematics, Columbia University, New York, NY 10027 \newline
\indent{\emailfont{ording@math.columbia.edu}}}
\thanks{Matthew Hedden was supported by an NSF postdoctoral fellowship and both authors received partial support from the NSF Holomorphic curves FRG grant during the course of this work}

\begin{document}

\begin{abstract}
In this paper we present several counterexamples to Rasmussen's conjecture that the concordance invariant coming from Khovanov homology is equal to twice the invariant coming from \os \ Floer homology.  The counterexamples are twisted Whitehead doubles of the $(2,2n+1)$ torus knots.
\end{abstract}

\include{intro}

\include{comp}

\commentable{
\bibliographystyle{plain}
\bibliography{biblio}
}

\end{document}

%% file: mac.tex
\newcommand\commentable[1]{#1}

\newtheorem{theorem}{Theorem}[section]
\newtheorem{prop}[theorem]{Proposition}

\newtheorem{lemma}[theorem]{Lemma}

\def\endproof{\relax\ifmmode\expandafter\endproofmath\else
  \unskip\nobreak\hfil\penalty50\hskip.75em\hbox{}\nobreak\hfil\bull
  {\parfillskip=0pt \finalhyphendemerits=0 \bigbreak}\fi}
\def\endproofmath$${\eqno\bull$$\bigbreak}
\def\bull{\vbox{\hrule\hbox{\vrule\kern3pt\vbox{\kern6pt}\kern3pt\vrule}\hrule}}

\newcommand{\Z}{\mathbb{Z}}

\newcommand{\OneHalf}{\frac{1}{2}}

\newcommand{\F}{\mathbb F}

\newcommand\relspinc{\underline{\spinc}}

\newcommand\ModSphere{\ModFlow\left({\mathbb S}\longrightarrow 
\Sym^{g-1}(\Sigma_{1})\times \Sym^2(\Sigma_{2})\right)}
\newcommand\ModSpheres\ModSphere

\newcommand\CFa{\widehat{CF}}

\newcommand\HFa{\widehat{HF}}

\newcommand\UnparModSp{\widehat \ModSp}
\newcommand\UnparModFlow\UnparModSp
\newcommand\Mod\ModSp

\newcommand{\spinc}{\mathfrak s}

\newcommand\ModMaps{\mathcal M}
\newcommand\ModSp\ModMaps

\newcommand\spincrel\relspinc

\newcommand\HFK{HFK}

\newcommand\HFKa{\widehat\HFK}

\newcommand\Dual{\mathcal D}
\newcommand\Duality\Dual

\newcommand\ons{Ozsv{\'a}th and Szab{\'o}}
\newcommand\os{Ozsv{\'a}th-Szab{\'o}}

%% file: intro.tex
\maketitle
\section{Introduction}

In \cite{FourBall} Ozsv{\'a}th and Szab{\'o} defined a smooth concordance invariant, denoted $\tau(K)$, whose value for the $(p,q)$ torus knot provided a new proof of Milnor's famous conjecture on the unknotting number of torus knots. Rasmussen independently discovered this invariant in his thesis, \cite{Ras1}. Milnor's conjecture has a long history in gauge theory, and its original proof is due to Kronheimer and Mrowka, \cite{KM}.  Recently, however, Rasmussen \cite{Ras2} discovered another smooth concordance invariant whose value for torus knots proves the conjecture.   Denoted $s(K)$, the invariant is defined using a refinement, due to Lee \cite{Lee}, of the purely combinatorial knot (co)homology  theory introduced by Khovanov \cite{Khovanov}.  Rasmussen's proof of the Milnor conjecture using $s$ is the first proof which avoids the analytical machinery of gauge theory.  It was noted immediately that the two invariants share several formal properties (e.g. an inequality relating the invariants of knots which differ by a crossing change) which in turn imply that they agree (or more precisely, that $s(K)$ and $2\tau(K)$ agree) for many knots.  For instance, $s(K)=2\tau(K)$ for the following families of knots:

\begin{enumerate}
	\item Torus knots:  $s(K)=2\tau(K)=2g(K)$ where $g(K)$ denotes the Seifert genus of $K$. This is due to Rasmussen \cite{Ras2}  for $s$ and \ons \ \cite{Lens} for $\tau$.  
	\item Alternating knots: $s(K)=2\tau(K)=\sigma(K)$ where $\sigma(K)$ is the classical Tristam-Levine signature of $K$. This is due to Lee \cite{Lee} for $s$, and \ons \ \cite{AltKnots} for $\tau$. 
	\item Strongly quasipositive knots, in particular positive knots: $s(K)=2\tau(K)=2g(K)$.  This is due to Livingston, \cite{Livingston1}.  See also \cite{Rudolph}.
	\item Quasipositive knots: $s(K)=2\tau(K)=2g_4(K)$, where $g_4(K)$ denotes the smooth slice genus of $K$. This follows from work of Plamenevskaya \cite{Olga1} for $\tau$ and from Plamenevskaya \cite{Olga2} and Shumakovitch \cite{Shumakovitch} for $s$. See also \cite{SQPfiber}.
	\item Knots with up to $10$ crossings.  \cite{Goda,Ras1,FourBall}
	\item ``Most'' twisted Whitehead doubles of an arbitrary knot, $K$.  This is due to Livingston and Naik \cite{Livingston2}
	\item Fibered knots with $\tau(K)=g(K)$. This follows from work of the first author \cite{SQPfiber}.
\end{enumerate}

\bigskip
\noindent Indeed, it was conjectured that the two invariants always coincide:

\bigskip
\noindent{\bf Conjecture}:({\em Rasmussen} \cite{Ras2}) $s(K)=2\tau(K)$ for all knots, $K$.  
\bigskip

\begin{figure}
\begin{center}

\caption{\label{fig:Double}
The $t$-twisted positive Whitehead double of the right-handed trefoil.  The box indicates the number of full right-handed twists to insert.}
\includegraphics{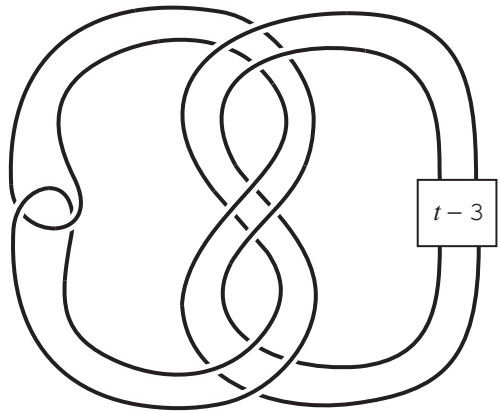}
\end{center}
\end{figure}

In light of the above list, the formal properties that the two invariants share, and several other striking connections between Khovanov's homology theory and  \os \ theory \cite{Seidel,Branched,Ciprian,Ras4}, there was justified hope that the above conjecture could be true.  However, we will demonstrate a counterexample:

\begin{theorem}
	Let $D_+(T_{2,3},2)$ denote the $2$-twisted positive Whitehead double of the right-handed trefoil knot (see Figure \ref{fig:Double}).  Then $\tau(D_+(T_{2,3},2)=0$ while $s(D_+(T_{2,3},2))=2$.  
\end{theorem}

Livingston and Naik \cite{Livingston2} calculate $\tau$ and $s$ for all but finitely many twisted Whitehead doubles of a knot, $K$, in terms of the maximal Thurston-Bennequin number of $K$, $TB(K)$, and its reflection, $\overline{K}$. In particular, they show that  $\tau(D_+(K,t))=s(D_+(K,t))/2=1$ if $t \le TB(K)$ and  $\tau(D_+(K,t))=s(D_+(K,t))/2=0$ if $t \ge -TB(\overline{K})$.  In light of an inequality satisfied by $\tau$ and $s$ under the operation of a crossing change, they define an invariant (which the results of this paper indicate is actually two invariants) $t_{\tau}(K)$ (resp. $t_s(K)$) which is the greatest integer $t$ such that $\tau(D_+(K,t))=1$ (resp. $s(D_+(K,t)=2$).  Using the techniques for the calculation above, we are able to determine $t_{\tau}(K)$ for the $(2,2n+1)$ torus knots:

\begin{theorem}
\label{thm:2n}
Let $D_+(T_{2,2n+1},t)$  denote the $t$-twisted positive Whitehead double of the the $(2,2n+1)$ knot.  Then we have:

$$ \tau(D_+(T_{2,2n+1},t))= 
\left\{\begin{array}{ll}

0  & {\text{for
 $t>2n-1$}} \\

 1 & {\text{for
 $t\le2n-1$}} \\

\end{array}
\right. $$
\end{theorem}

Thus, $t_{\tau}(T_{2,2n+1})=2n-1$.  In fact, the above knots provide further counterexamples, as was shared with us by Jake Rasmussen, who used Bar-Natan's  program \cite{BarNatan} for computing Khovanov homology to calculate $s$ for the knots in the above family which are not covered by Livingston and Naik's result.  In particular:
$$ s(D_+(T_{2,5},5))= s(D_+(T_{2,5},4))= s(D_+(T_{2,7},8))= s(D_+(T_{2,7},7))= s(D_+(T_{2,7},6))=2, $$
\noindent while Theorem \ref{thm:2n} implies that $\tau=0$ for these knots.  It seems likely that Whitehead doubles of the $(2,2n+1)$ torus knots provide an infinite family of counterexamples.  Indeed, it would be reasonable to guess that $t_s(T_{2,2n+1})=3n-1$

We prove the above results first by calculating the knot Floer homology groups of a specific twisted Whitehead double which happens to be a $(1,1)$ knot. A general technique for calculating the Floer homology of such knots was developed by Goda, Morifuji, and Matsuda \cite{Goda} and we apply their technique here.  We then use results of Eftekhary \cite{Eaman} for the $0$-twisted Whitehead double of  $T_{2,2n+1}$, together with properties of the skein exact sequence for knot Floer homology to calculate $\tau$ for the examples above.  The techniques here will be refined and generalized in \cite{Doubling} to calculate $\tau$ and some of the Floer homology of an arbitrarily twisted Whitehead double of an arbitrary knot (in fact, \cite{Doubling} will prove that $t_{\tau}(K)=2\tau(K)-1$).  We also remark that $(1,1)$ satellite knots were classified by Morimoto and Sakuma in \cite{Morimoto}, and it was in the context of a more general study of these knots that this work arose.   We hope to return to this study (see also \cite{PhilThesis}).

	We conclude by remarking that there is a beautiful conjectural picture  due to Dunfield, Gukov, and Rasmussen, \cite{Dunfield} of a triply graded homology theory which would unify Khovanov homology, knot Floer homology, and the various $sl(n)$ link homology theories of Khovanov and Rozansky \cite{Khovanov2}.  It would be very interesting to understand this conjecture for the above examples - in particular it would be useful to calculate the $sl(n)$ link homology.
	
	\noindent{\bf{Acknowledgements:}} It is our pleasure to thank  Hiroshi Goda, Hiroshi Matsuda, Peter Ozsv{\'a}th, Zoltan Szab{\'o}, and Jake Rasmussen for many stimulating discussions.  We owe special thanks to Jake for computing $s$ for the examples in this paper and for sharing both his results and his general knowledge of the $s$ invariant with us.

%% file: comp.tex
\section{Computation}
\label{sec:comp}

\subsection{Computation of s}
We begin by computing $s$.  Using Bar-Natan and Shumakovitch's programs \cite{BarNatan,Shumakovitch2} to compute Khovanov homology, we obtained the following Poincar{\'e} polynomial for the unreduced Khovanov homology of $D_+(T_{2,3},2)$:

\begin{eqnarray*}
PKh(q,t)=q^{-5}t^{-4}+q^{-1}t^{-3}+q^{-1}t^{-2}+qt^{-1}+q^3t^{-1}+2q+q^3+q^5+2q^5t+
q^{5}t^2+
q^{9}t^{2}+\\
q^{7}t^{3}+
q^{9}t^{3}+
q^{7}t^{4}+
q^{11}t^{4}+
q^{9}t^{5}+
q^{11}t^{5}+
q^{13}t^{6}+
q^{13}t^{7}+
q^{15}t^{8}+
q^{17}t^{8}+
q^{19}t^{9}.
\end{eqnarray*}

The only homology in homological grading $0$ is supported in $q$ gradings $1,3,5$.  It follows from the definition of $s$ that $s(D_+(T_{2,3},2))$ is equal to $2$ or $4$.  However, the fact that the genus of  $D_+(T_{2,3},2)$ is equal to one and $|s(K)|\le 2g_4(K)$ implies  $s(D_+(T_{2,3},2))=2$.
\subsection{Computation of $\tau$}
We begin our calculations by computing the knot Floer homology groups of $D_+(T_{2,3},6)$, the $6$-twisted positive Whitehead double of the right-handed trefoil.  For notational simplicity, we hereafter denote the $t$-twisted positive Whitehead double of the right-handed trefoil by $D(t)$.  

\begin{prop}
\label{prop:D6}
$$ \HFKa_*(D(6),i)\cong 
\left\{\begin{array}{ll}

	\F_{(1)}^{4}\oplus\F_{(-1)}^2  & {\text{for
 $i=1$}}\\

 \F_{(0)}^{9}\oplus\F_{(-2)}^4  & {\text{for
 $i=0$}}\\

	\F_{(-1)}^{4}\oplus\F_{(-3)}^2  & {\text{for
 $i=-1$}}\\

\end{array}
\right. $$
Where $\F$ denotes the field with 2 elements.

\end{prop}

\noindent{\bf Remark:} Note that $\tau(D(6))=0$. There is simply no homology in grading $0$ supported in filtration grading $1$ or $-1$.

\begin{proof}
	We first apply the technique developed in \cite{Goda} for obtaining a genus one doubly-pointed Heegaard diagram from a $(1,1)$ presentation to the knot at hand, $D(6)$.  This is illustrated in Figure \ref{fig:Unwind}.  Following the technique which \ons \ introduced in Section $6$ of \cite{Knots} (and which was further developed by \cite{Goda}), we lift this genus one diagram to the universal cover, Figure \ref{fig:Lift}, and compute the boundary map:

\begin{eqnarray*}
\partial[x_1,i,i]&=&0 \\
\partial[x_2,i,i+1]&=&[x_1,i,i]+[x_5,i-1,i-1] \\
\partial[x_3,i,i]&=& [x_2,i-1,i]+[x_4,i,i-1] \\
\partial[x_4,i,i-1]&=& [x_1,i-1,i-1]+[x_5,i-2,i-2] \\
\partial[x_5,i,i]&=& 0 \\
\partial[x_6,i,i+1]&=&[x_5,i,i] +[x_9,i,i]\\
\partial[x_7,i,i]&=& [x_6,i-1,i]+[x_8,i,i-1] \\
\partial[x_8,i,i-1]&=&[x_5,i-1,i-1]+[x_9,i-1,i-1] \\
\partial[x_9,i,i] &=& 0 \\
\partial[x_{10},i,i+1]&=&[x_9,i,i] +[x_{13},i,i]\\
\partial[x_{11},i,i]&=& [x_{10},i-1,i]+[x_{12},i,i-1] \\
\partial[x_{12},i,i-1]&=&[x_9,i-1,i-1]+[x_{13},i-1,i-1] \\
\partial[x_{13},i,i] &=& 0 \\
\partial[x_{14},i,i+1]&=&[x_{13},i,i] +[x_{17},i,i]\\
\partial[x_{15},i,i]&=& [x_{14},i-1,i]+[x_{16},i,i-1] \\
\partial[x_{16},i,i-1]&=&[x_{13},i-1,i-1]+[x_{17},i-1,i-1] \\
\partial[x_{17},i,i] &=& 0 \\
\partial[x_{18},i,i+1]&=&[x_{17},i,i] +[x_{21},i,i]\\
\partial[x_{19},i,i]&=& [x_{18},i-1,i]+[x_{20},i,i-1] \\
\partial[x_{20},i,i-1]&=&[x_{17},i-1,i-1]+[x_{21},i-1,i-1] \\
\partial[x_{21},i,i]&=&0 \\
\partial[x_{22},i,i+1]&=&[x_{25},i,i]+[x_{21},i-1,i-1] \\
\partial[x_{23},i,i]&=& [x_{22},i-1,i]+[x_{24},i,i-1] \\
\partial[x_{24},i,i-1]&=& [x_{25},i-1,i-1]+[x_{21},i-2,i-2] \\
\partial[x_{25},i,i]&=& 0. \\
\end{eqnarray*}

\begin{figure}
\begin{center}
\caption{\label{fig:Unwind}
Illustration of the process by which we obtain a doubly-pointed Heegaard diagram for $D(6)$, Step (f), from its $(1,1)$ presentation, Step (a).  }
\includegraphics{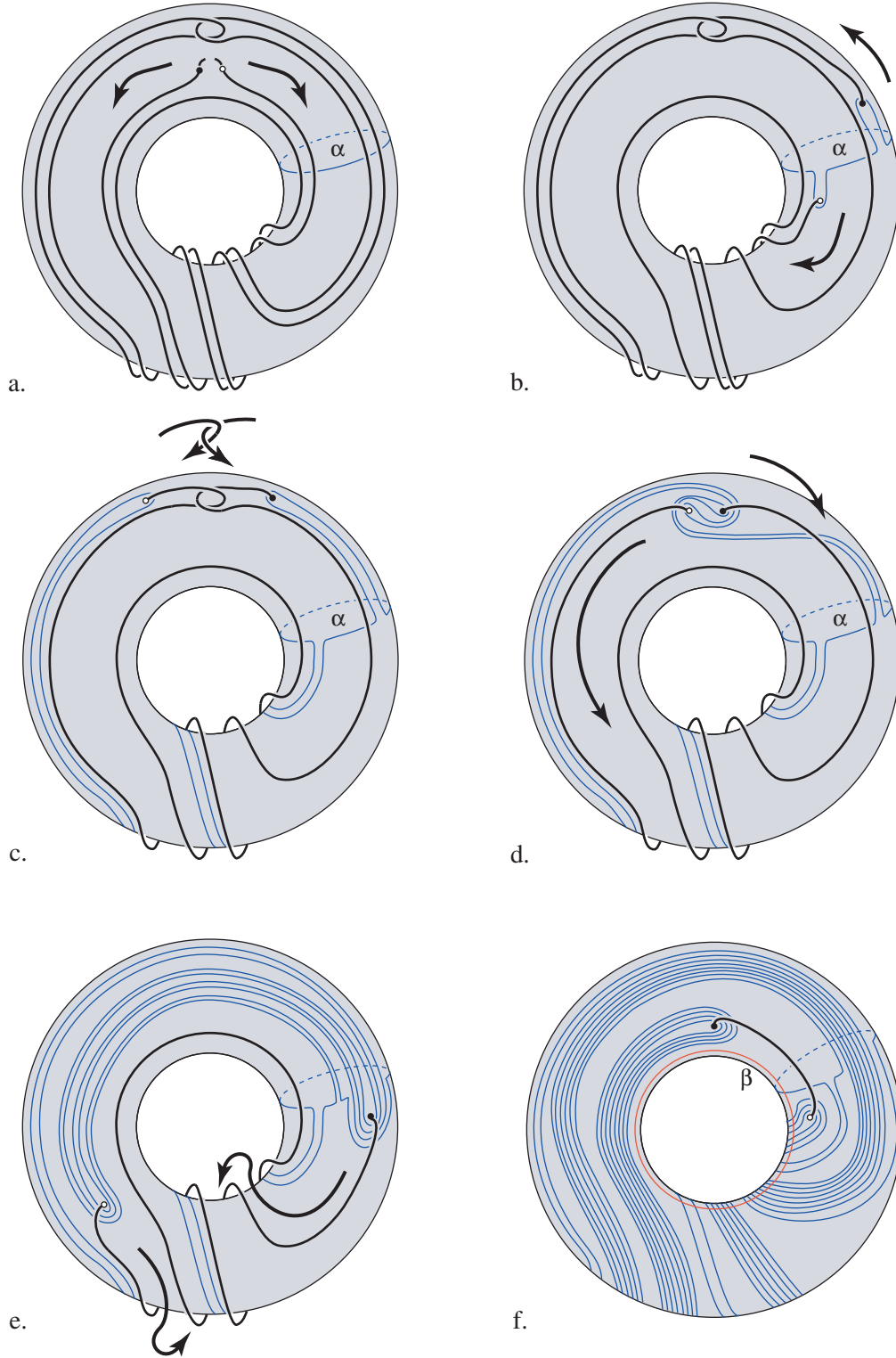}
\end{center}
\end{figure}

\begin{figure}
\begin{center}

\caption{\label{fig:Lift}
The Heegaard diagram of the previous figure, lifted to the universal cover of the torus.  We have chosen a particular lift of $\alpha$ and $\beta$, as indicated.  The open circles denote lifts of the basepoint $z$ while the black circles denote lifts of $w$. \bigskip}
\includegraphics{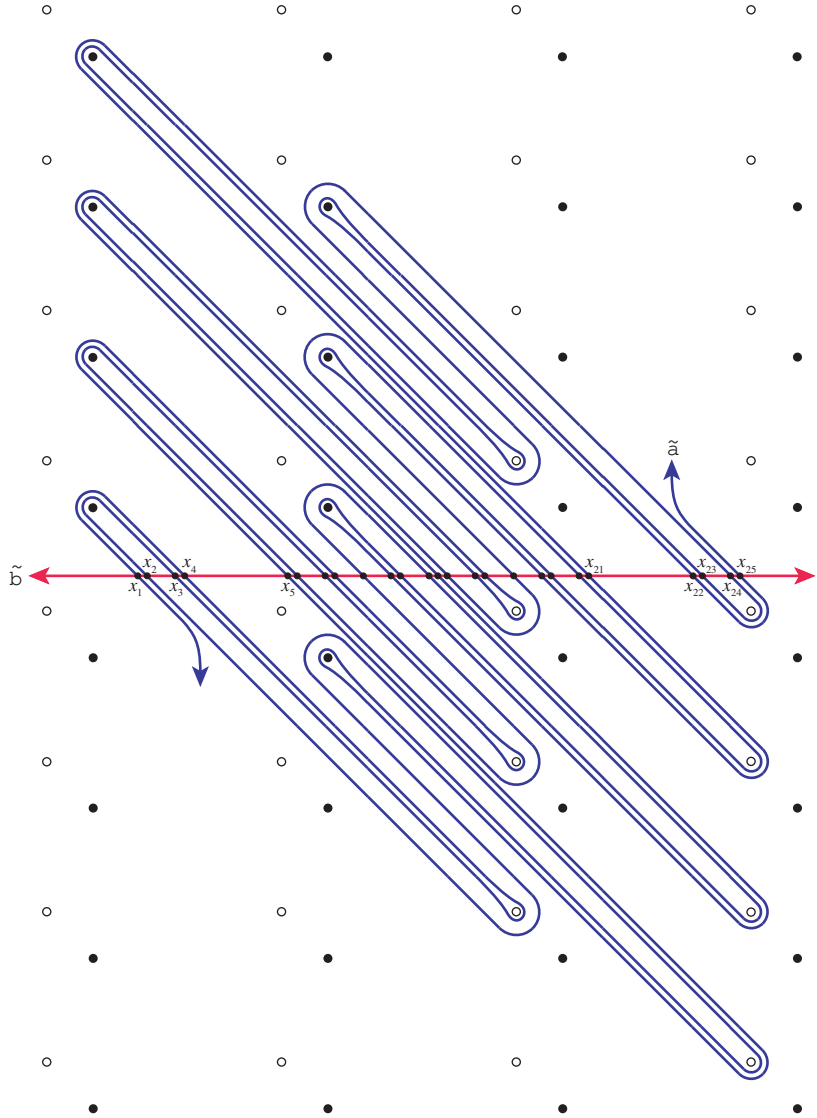}
\end{center}
\end{figure}

\newpage
\noindent Using our knowledge of the differential, it is easy to separate the generators of the chain complex into their respective filtration and homological gradings.  In the table below the vertical (horizontal) direction indicates the filtration (homological) grading:

\begin{center}
\begin{tabular}{|c||c|c|c|c|c|}
\noalign{\hrule height 1.1pt}
 & -3 & -2 & -1 & 0 & 1  \\

\hline
$1$ & & & $x_2,x_{22}$ & & $x_6,x_{10},$ \\
&&& && $x_{14},x_{18}$ \\
\hline
$0$ & &$x_1,x_3,$ & & $x_5,x_7,x_9,x_{11},x_{13},$ &  \\
 & & $x_{23},x_{25}$& &  $x_{15},x_{17},x_{19},x_{21}$ & \\
\hline
$-1$ & $x_4,x_{24}$ & & $x_8,x_{12},$ & & \\
&&& $x_{16},x_{20}$ && \\
\noalign{\hrule height 1.1pt}
\end{tabular}
\end{center}
\bigskip
The proposition follows immediately.
\end{proof}

\bigskip
\noindent Next we recall the following result of Eftekhary \cite{Eaman}:

\begin{theorem} 
	\label{thm:D0} (Eftekhary \cite{Eaman})
	$$\HFKa_*(D(0),1) \cong  \F^2_{(m)}\oplus \F^2_{(m-1)},$$ where the subscript $(m)$ indicates that the homological degree is known only as a relative $\Z$-grading.
\end{theorem}

By performing $6$ successive crossing changes to the twisting region of the knot diagram shown in Figure $1$, we can change $D(6)$ into $D(0)$.  Each of these operations changes a negative crossing to a positive crossing.   There is a skein exact sequence for each crossing change which relates the Floer homology groups of $D(t)$, $D(t-1)$ and the two-component link obtained from the oriented resolution of the crossing which we change.  For each $t$, this latter link is the positive Hopf link, which we denote by $H$.  The Floer homology of $H$ (i.e. the Floer homology of its ``knotification'', see Section $2$ of \cite{Knots}) is given by:

\begin{prop}
\label{prop:H}

$$\HFKa(H,i)\cong
\left\{\begin{array}{ll}

\F_{(\OneHalf)} & {\text{if $i=1$}} \\
\F^2_{(-\OneHalf)} & {\text{if $i=0$}} \\
\F_{(-\frac{3}{2})} & {\text{if $i=-1$}} \\
0 & {\text{otherwise}} \\
\end{array}
\right. $$

\end{prop}
\begin{proof} This was originally proved in Proposition $9.2$ of \cite{Knots}, but also follows easily from \cite{AltKnots}, whose main theorem determines the Floer homology of alternating links in terms of their Alexander polynomial and signature.
\end{proof}

The remaining step in our computation of $\tau(D(2))$ will be to study the skein exact sequences corresponding to the $6$ aforementioned crossing changes.  In each case, the skein sequence for the top filtration level takes the following form:

$$
\begin{CD}
	... @>>>	\HFKa(D(t),1)
	@>{f_1}>> \F_{(\OneHalf)} @>{f_2}>> 	 \HFKa(D(t-1),1)  @>{f_3}>>...
\end{CD}
$$

Where the maps $f_1$ and $f_2$ lower homological degree by one-half and $f_3$ is non-increasing in the homological degree.  We also note that \ons \ define an absolute $\Z/2\Z$ grading on the groups in the above sequence which is simply the parity of the $\Z$-grading, with the convention that the $\F_{(\OneHalf)}$-summand is supported in odd parity.  With respect to the $\Z/2\Z$ grading, the maps $f_1$ and $f_3$ are grading-preserving, while $f_2$ is grading-reversing.  It follows at once that there are two options for each skein sequence:

\begin{enumerate}
	\item $f_2$ is trivial, $f_1$ is non-trivial
	\item $f_2$ is non-trivial, $f_1$ is trivial
\end{enumerate}
We make the following claim:

\begin{prop}
	\label{prop:skein}
 In the exact sequence above relating $D(t)$, $D(t-1)$ and $H$, the map $f_2$ is non-trivial if and only if  $\tau(D(t-1))=1$.  Furthermore, if $\tau(D(t-1)\ne1$, then it is equal to $0$.
\end{prop}

\noindent {\bf Remark:} The second part of the statement also follows from work of Livingston and Naik \cite{Livingston2}.

\begin{proof} This will follow from the fact that $f_2$ is the lowest order term in a filtered chain map, $\tilde{f_2}$, between chain complexes which are chain homotopy equivalent to $\CFa(S^1\times S^2)$ and $\CFa(S^3)$, respectively.
	
To begin, note that the Floer homology groups for $H$ (resp. $D(t-1)$) are endowed with an induced differential which gives them the structure of a filtered chain complex. Moreover, this differential strictly lowers the filtration index. In the case of $H$, the homology of this filtered chain complex is  $\HFa(S^1\times S^2)\cong \F_{(-\OneHalf)}\oplus \F_{(\OneHalf)}$.  In the case of $D(t-1)$, the homology is $\HFa(S^3)\cong \F_{(0)}$.  The filtration on the knot Floer homology of $D(t-1)$ induces a filtration on $\HFa(S^3)$ in the standard way i.e. the filtration level of any cycle, $z = \Sigma n_x x$, is by definition the maximum filtration level of any chain $x$ which comprises $z$.   Now $\tau(D(t-1))$ is defined to be the minimum filtration degree of any cycle $z\in \HFKa(D(t-1))$ which is homologous to a generator of $\HFa(S^3)$.  

It follows from the proof of the skein sequence (Theorem $8.2$ of \cite{Knots}) that there is a map $$\tilde{f_2} :\HFKa(H) \rightarrow \HFKa(D(t-1)),$$

\noindent which commutes with the differentials on $\HFKa$ and respects the filtration i.e. does not increase the filtration index. In other words, $\tilde{f_2}$ is a filtered chain map between filtered chain complexes.  Furthermore, $\tilde{f_2}$ decomposes as a sum of homogeneous pieces, each of which lower the filtration by some fixed integer.  The map in the skein sequence is the part of $\tilde{f_2}$ which preserves (does not lower) the filtration, restricted to each filtration, $\HFKa(H,i)$ (in the case at hand $i=1$).  

From Proposition \ref{prop:H} we see that a chain generating $\HFKa(H,1)\cong \F_{(\OneHalf)}$  is a cycle under the induced differential, and hence the above discussion implies that $\tilde{f_2}$  maps this chain to a cycle,  $z\in \HFKa(D(t-1))$.  Now if $f_2$ is non-trivial, $z$ contains non-trivial chains with filtration index $1$.  The definition of $\tau$, together with the fact that $\HFKa(D(t-1),i)\cong 0$ for $i>1$ implies $\tau(D(t-1))=1$. 

Now on the level of homology, $\tilde{f_2}$ induces a map:  
$$
\begin{CD}
	\HFa(S^1\times S^2)\cong \F_{(-\OneHalf)}\oplus \F_{(\OneHalf)} @>{(\tilde{f_2})_*}>> 	\HFa(S^3)\cong \F_{(0)}  \\
\end{CD}
$$
\noindent which sends the space supported in degree one-half to the generator.  If $\tau(D(t-1))=1$, the cycle  generating $\HFa(S^3)$  contains non-trivial chains in filtration level $1$.  It follows that $f_2$ - the part of $\tilde{f_2}$ which preserves the filtration - is non-trivial.  

 Finally, if $\tau(D(t))=-1$, a similar analysis shows that $\tilde{f_1}$ restricted to $\HFKa(D(t),-1)$ would raise the filtration degree, contradicting the fact that this map respects the filtration.

\end{proof}

The above proposition shows that the map $f_2$ in the skein sequence controls the behavior of $\tau(D(t-1))$.  We determine when $f_2$ is non-trivial in the six applications of the sequence:  

\begin{lemma}
	\label{lemma:lemma1}
The map	$ f_2:\HFKa(H,1)\rightarrow \HFKa(D(t-1),1)$ is trivial for $t=6,5,4,3$ and non-trivial for $t=2,1$.
\end{lemma}

The theorem about $\tau(D(2))$ will follow immediately from the above lemma and Proposition \ref{prop:skein}.  Indeed, it follows easily from the proof that $\tau(D(t))=0$ if $t>1$ and  $\tau(D(t))=1$ if $t\le 1$.  

\begin{proof}  We study following function: $$e(t)=\mathrm{rk}_{even} \HFKa(D(t),1),$$

\noindent which measures the rank of the Floer homology in top filtration level supported in even homological degree.	

\noindent {\bf Claim:}  If $f_2$ is non-trivial then $e(t-1)=e(t)+1$. If $f_2$ is trivial then $e(t-1)=e(t)$.

The claim follows from the form of the skein sequence at hand, together with the knowledge that $f_1$ and $f_3$ preserve the $\Z/2\Z$-grading while $f_2$ reverses it.  It follows from Proposition \ref{prop:D6} and Theorem \ref{thm:D0} that $e(6)=0$ and $e(0)=2$.  Thus the claim shows that among the six applications of the skein sequence, $f_2$ is non-trivial exactly twice.  
	
	Next, recall that $\tau$ (and $s$) satisfy the following inequality under the operation of changing a crossing in a given knot diagram (see \cite{Livingston1} or \cite{FourBall} for a proof):

$$\tau(K_+)-1\le \tau(K_-) \le \tau(K_+),$$

\noindent where $K_+$ (resp. $K_-$) denote the diagram with the positive (resp. negative) crossing. 
Now each application of the skein sequence arose from changing a single negative crossing to a positive crossing.  Hence the above inequality becomes (for $k>0$):

$$\tau(D(t-k))-k\le \tau(D(t)) \le \tau(D(t-k)).$$

If $f_2$ were non-trivial for some $t$ and trivial for $t-k$, then Proposition \ref{prop:skein} would imply $\tau(D(t-1))=1$ and $\tau(D(t-k-1)= 0$ violating the inequality.  Thus $f_2$ is trivial for $t=6,5,4,3$ as stated, and non-trivial for $t=2,1$.  
\end{proof}

\subsection{Twisted Whitehead doubles of $(2,2n+1)$ torus knots}

Let $D_+(T_{2,2n+1},t)$  denote the $t$-twisted positive Whitehead double of the right-handed $(2,2n+1)$ torus knot. Results of \cite{Morimoto} indicate that the $D_+(T_{2,2n+1},4n+2)$ is a $(1,1)$ knot, and indeed we can repeat the calculation of Proposition \ref{prop:D6} to yield:






\begin{prop}
\label{prop:Kstart}
$$ \HFKa_*(D_+(T_{2,2n+1},4n+2),i)\cong 
\left\{\begin{array}{ll}

	\F_{(1)}^{2n+2}\oplus\F_{(-1)}^2\oplus\F_{(-3)}^2\ldots \oplus\F_{(-2n+1)}^2  & {\text{for
 $i=1$}} \\

 \F_{(0)}^{4n+5}\oplus\F_{(-2)}^4\oplus\F_{(-4)}^4\ldots \oplus\F_{(-2n)}^4  & {\text{for
 $i=0$}} \\

	\F_{(-1)}^{2n+2}\oplus\F_{(-3)}^2\oplus\F_{(-5)}^2\ldots \oplus\F_{(-2n-1)}^2  & {\text{for
 $i=-1$}} \\

\end{array}
\right. $$
\end{prop}

In addition, Eftekhary's \cite{Eaman} results in this case yield:

\begin{theorem}
	\label{prop:eaman}
	
$$ \HFKa_*(D_+(T_{2,2n+1},0),1)\cong \F_{(m)}^{2n}\oplus\F_{(m-1)}^2\oplus\F_{(m-3)}^2\ldots \oplus\F_{(m-2n+1)}^2 $$

\end{theorem}

The technique for computing $\tau$ in the case of the trefoil can now be applied to yield Theorem \ref{thm:2n}.  This result should be compared with results of Livingston and Naik \cite{Livingston2}.